\documentclass[12 pt, a4 paper]{article}
\usepackage{graphicx}
\usepackage{epstopdf}
\usepackage{amsfonts}
\usepackage{amssymb}
\usepackage{amsfonts,amsmath,latexsym,amssymb}                                                                      
\usepackage{amsthm}
\usepackage{cite}
\usepackage{lipsum} 
\usepackage{float} 
\usepackage{hyperref}
\hypersetup{
colorlinks,
citecolor=red,
linkcolor=blue,
urlcolor=pink}
\usepackage{mathrsfs,upref}
 
\numberwithin{equation}{section}
\newtheorem{theorem}{Theorem}

\newtheorem{lemma}{Lemma}

\DeclareMathOperator{\tanhinv}{\tanh^{-1}}
\begin{document}

\begin{center}\begin{Large} 
\textbf{Alternative proofs of Shafer's inequality for inverse hyperbolic tangent}
\end{Large}\end{center}
\begin{center}
Yogesh J. Bagul$^{1,*}$,  Ramkrishna M. Dhaigude$^{2}$  \\
$^{1,*}$Department of Mathematics, \\
K. K. M. College Manwath,\\
Dist: Parbhani(M.S.) - 431505, India.\\
Email: yjbagul@gmail.com\\
$^{2}$Department of Mathematics, \\
Government Vidarbha Institute of Science,\\
and Humanities, Amravati(M. S.)-444604, India\\
Email: rmdhaigude@gmail.com
\end{center} 

\textbf{Abstract.}
We point out that a concise proof of Theorem 2 in the article, 
{\it 'On a quadratic estimate of Shafer'} by L. Zhu contains a small mistake. Correcting this mistake and giving alternative proofs of Theorem 2 is the main aim of this note. \\

\renewcommand{\thefootnote}{}
\footnotetext{ $^*$Corresponding author
\par

}


\section{Introduction and Correction}

In 2008, L. Zhu \cite{zhu} published a new proof of the following theorem:

\begin{theorem}\label{thm1}
Let $ 0 < x < \sqrt{15}/4. $ Then
\begin{align}\label{eqn1.1}
\frac{\tanhinv x}{x} < \frac{8}{3 + \sqrt{25 - \frac{80}{3}x^2}}.
\end{align}
\end{theorem} 
The inequality (\ref{eqn1.1}) was originally established by R. E. Shafer \cite{shafer, shafer1, shafer2} and its alternative proof is given in \cite{zhu} in a concise way. Though the proof of Theorem \ref{thm1} is given  in a simple way in \cite{zhu}, it contains a small mistake which can be explained as follows:

 While giving the proof of Theorem \ref{thm1}, it is shown in \cite{zhu}  that the function
$$ H(x) = \frac{25-\left(\frac{8x}{\tanhinv x} - 3\right)^2}{x^2} $$
is decreasing on $(0, \sqrt{15}/4). $ This is accomplished by showing 
$$ I(t) = \frac{-4\sinh^2 t + 3 t \sinh t \cosh t + t^2 \cosh^2 t}{t^4 \cosh^2 t} = \frac{A(t)}{B(t)} $$
to be decreasing on $ (0, \tanhinv \sqrt{15}/4) $ due to the transformation $ H(x) = 16 I(t),$ where $ \tanhinv x = t. $ A careful observation shows that the denominator $ B(t) $ of $ I(t) $ is mistaken as $ t^4 \cosh^2 t $ instead of $ t^2 \sinh^2 t. $ Fortunately, the function $ I(t) $ remains decreasing for either expression for $ B(t)$ and the final conclusion is unaffected. For final conclusion, the following lemma is used.

\begin{lemma}(\cite{heikkala})  \label{lemm1} 
Let $ A(x) = \sum_{n = 0}^{\infty}a_n x^n $ and $ B(x) = \sum_{n = 0}^{\infty}b_n x^n $ be convergent for $ \vert x \vert < R ,$ where $ a_n $ and $ b_n $ are real numbers for $ n = 0, 1, 2, \cdots $ such that $ b_n > 0 .$ If the sequence ${a_n/b_n}$ is strictly increasing(or decreasing), then the function $ A(x)/B(x) $ is also strictly increasing(or decreasing) on $ (0, R).$ 
\end{lemma}

However, it is necessary  to show that how $ I(t) $ is decreasing on $ (0, \tanhinv \sqrt{15}/4) $ with $ B(t) = t^2 \sinh^2 t. $ In fact, with this $ B(t) $ the proof becomes more clear and convincing. Here we present the proof. 

{\it Corrected proof of Theorem \ref{thm1}.} As in the concise proof of Theorem 2 in \cite{zhu}, we have 
$$ A(t) = \sum_{n=1}^{\infty}a_n t^{2n+2}, $$ where
$$ a_n = \frac{-4\cdot2^{2n+2} + 3(2n +2)\cdot 2^{2n+1} + (2n+1)(2n+2)\cdot 2^{2n}}{2(2n+2)!}. $$
 Now 
\begin{align*}
B(t) &= t^2 \sinh^2 t \\
&= \frac{t^2}{2}\left(\cosh 2t - 1\right) \\
&= - \frac{1}{2} t^2 + \frac{1}{2} t^2 \sum_{n=0}^{\infty}\frac{2^{2n}}{(2n)!}t^{2n} \\
&= \sum_{n=1}^{\infty}\frac{2^{2n-1}}{(2n)!}t^{2n+2} = \sum_{n=1}^{\infty}b_n t^{2n+2}
\end{align*}
where $ b_n = \frac{2^{2n-1}}{(2n)!} = \frac{(n+1)(2n+1)\cdot2^{2n}}{(2n+2)!}. $ Then we write 
\begin{align*}
\frac{a_n}{b_n} &= \frac{(n+1)(2n+1)+6(n+1)-8}{(n+1)(2n+1)} \\
&= \frac{2n^2 + 9n -1}{2n^2 + 3n +1} \\
&= 1 + \frac{2(3n-1)}{2n^2+3n+1} := 1 + 2 c_n
\end{align*} 
where $$ c_n = \frac{3n - 1}{2n^2 + 3n +1} \, \,  \text{and} \, \,  c_{n+1} =\frac{3n + 2}{2n^2 + 7n + 6}, \, \, n = 1, 2, 3, \cdots $$
We claim that 
$$ c_n \geq c_{n+1} , \, \, n \geq 1. $$
Equivalently, $$ \frac{3n - 1}{2n^2 + 3n + 1} \geq \frac{3n + 2}{2n^2 + 7n +6} $$ or $$ 19n^2 + 11n -6 \geq 13n^2 + 9n + 2. $$ i.e., $ 6n^2 + 2n \geq 8 $ which is true for $ n \geq 1. $ Therefore a sequence $ \left\lbrace \frac{a_n}{b_n} \right\rbrace $ is decreasing for $ n = 1, 2, 3, \cdots. $ Hence by  Lemma \ref{lemm1}, $ I(t) $ is also decreasing on 
$ (0, \tanhinv \sqrt{15}/4). $ $ \qed $ \\

Next, it is interesting to see other simple proofs of Theorem \ref{thm1}.
\section{Alternative simple proofs}

We give two alternative simple proofs of Theorem \ref{thm1}. The first proof, is very elementary and uses basic calculus only.

{\bf First simple proof of Theorem \ref{thm1}.} If we let $ \tanhinv x = t, $ then it suffices to prove that 
$$ \frac{t}{\tanh t} < \frac{8}{3 + \sqrt{25 - \frac{80}{3}\tanh^2 t}}, $$ for $ t \in (0, \tanhinv \sqrt{15}/4).$ Equivalently we want
$$\left( 8 \frac{\sinh t}{t} - 3 \cosh t \right)^2 > 25 \cosh^2 t - \frac{80}{3}\sinh^2 t. $$ i.e.
$$  64 \sinh^2 t - 48 t \sinh t \cosh t > 16 t^2 \cosh^2 t - \frac{80}{3} t^2 \sinh^2 t.  $$ Or
$$  192 \sinh^2 t - 144 t \sinh t \cosh t > 48 t^2 \cosh^2 t - 80 t^2 \sinh^2 t .  $$ i.e.
$$ 12 \sinh^2 t - 9t \sinh t \cosh t > 3t^2 - 2t^2 \sinh^2 t.  $$ 
Now suppose, $$ f(t) = 12 \sinh^2 t + 2t^2 \sinh^2 t - 9t \sinh t \cosh t - 3t^2. $$ Successive differentiations with respect to $ t $ give
\begin{align*}
f'(t) &= 15 \sinh t \cosh t - 5 t \sinh^2 t + 4t^2 \sinh t \cosh t - 9t \cosh^2 t - 6t, \\
f''(t) &= 10 \sinh^2 t + 15 \cosh^2 t - 20t \sinh t \cosh t + 4t^2 \sinh^2 t + 4 t^2 \cosh^2 t \\
&- 9 \cosh^2 t - 6, \\
f'''(t) &= 12 \sinh t \cosh t - 12 t \sinh^2 t - 12 t \cosh^2 t + 16 t^2 \sinh t \cosh t, \\
f^{iv}(t) &= 16 t^2 \sinh^2 t + 16 t^2 \cosh^2 t - 16 t \sinh t \cosh t \\
&= 16 t^2 \sinh^2 t + 16 t \cosh t ( t \cosh t - \sinh t) > 0
\end{align*}
due to well-known inequality $ \frac{\sinh t}{t} < \cosh t, \, \, t > 0. $ This implies that $ f'''(t) $ is strictly increasing for $ t > 0 $ and hence $ f'''(t) > f'''(0). $ Since, $ f'''(0) = f''(0) = f'(0) = f(0),$  we continue the argument and conclude that $ f(t) > f(0) = 0. $ This completes the proof. \\

{\bf Second simple proof of Theorem \ref{thm1}.} Since $$ I(t)   = \frac{t^2 \cosh^2 t + 3 t \sinh t \cosh t -4\sinh^2 t}{t^2 \sinh^2 t}  $$ we have 
 $$ I'(t) =  -\frac{1}{4t^3 \sinh^3 t} i(t), $$
where $$ i(t) = \left(24 \sinh t - 8 \sinh 3t - 3t \cosh t + 3t \cosh 3t + 8t^3 \cosh t + 12t^2 \sinh t\right). $$
Substituting the two formulas $$ \cosh kt = \sum_{n=0}^{\infty} \frac{k^{2n}}{(2n)!} t^{2n} \, \, \text{and} \, \, \sinh kt = \sum_{n=0}^{\infty} \frac{k^{2n+1}}{(2n+1)!} t^{2n+1}$$ into the previous formula to get

\begin{align*}
i(t) &= 24 \sum_{n=0}^{\infty}\frac{1}{(2n+1)!}t^{2n+1} - 8\sum_{n=0}^{\infty}\frac{3^{2n+1}}{(2n+1)!}t^{2n+1} - 3t\sum_{n=0}^{\infty}\frac{1}{(2n)!}t^{2n} \\
 &+ 3t\sum_{n=0}^{\infty}\frac{3^{2n}}{(2n)!}t^{2n} + 8t^3 \sum_{n=0}^{\infty}\frac{1}{(2n)!}t^{2n} + 12t^2 \sum_{n=0}^{\infty}\frac{1}{(2n+1)!}t^{2n+1} \\
 &= 24 \sum_{n=0}^{\infty}\frac{1}{(2n+1)!}t^{2n+1} - 8\sum_{n=0}^{\infty}\frac{3^{2n+1}}{(2n+1)!}t^{2n+1} - 3\sum_{n=0}^{\infty}\frac{1}{(2n)!}t^{2n+1} \\
 &+ \sum_{n=0}^{\infty}\frac{3^{2n+1}}{(2n)!}t^{2n+1} + 8 \sum_{n=0}^{\infty}\frac{1}{(2n)!}t^{2n+3} + 12 \sum_{n=0}^{\infty}\frac{1}{(2n+1)!}t^{2n+3} \\
 &=  \sum_{n=1}^{\infty}\frac{24}{(2n+1)!}t^{2n+1} - \sum_{n=1}^{\infty}\frac{8 \times 3^{2n+1}}{(2n+1)!}t^{2n+1} - \sum_{n=1}^{\infty}\frac{3}{(2n)!}t^{2n+1} \\
 &+ \sum_{n=1}^{\infty}\frac{3^{2n+1}}{(2n)!}t^{2n+1} +  \sum_{n=1}^{\infty}\frac{8}{(2n-2)!}t^{2n+1} +  \sum_{n=1}^{\infty}\frac{12}{(2n-1)!}t^{2n+1} \\
 &= \sum_{n=4}^{\infty}\frac{d_n}{(2n+1)!}t^{2n+1}
\end{align*}
where 
\begin{align*}
d_n &= 24 - 8 \times 3^{2n+1} - 3(2n+1) + 3^{2n+1}(2n+1) + 8(2n-1)(2n)(2n+1) \\
 &+ 12(2n)(2n+1) \\
&= 24 - 8 \times 3^{2n+1} -6n -3 + 3^{2n+1} \cdot 2n + 3^{2n +1} +16n(4n^2 - 1) \\
 &+ 24n (2n +1) \\
 &= 3^{2n+1} \cdot (2n-7) + 64n^3 + 48n^2 + 2n +21 > 0.
\end{align*}
So $ i(t) > 0 $ holds for all $ t > 0 $ giving us $ I'(t) < 0.$ Thus $ I(t) $ is decreasing on $ (0, \tanhinv \sqrt{15}/4) $ and so is $ H(x) $ on $ (0, \sqrt{15}/4). $ Consequently, 
$$ H(0+) > H(x) $$ and with $ H(0+) = 80/3 $  we get the inequality (\ref{eqn1.1}). $ \qed $ \\

{\bf Competing Interests.}
Authors would like to state that they do not have any competing interest.

{\bf Author's contributions.} Both the authors contributed equally for this paper.

{\bf Acknowledgement.} The authors would like to sincerely thank Prof. L. Zhu for checking the manuscript before submission and suggesting second simple proof.

\end{document}